\documentclass[10pt]{article}

\newtheorem{theorem}{Theorem}
\newtheorem{proposition}{Proposition}
\newtheorem{lemma}{Lemma}

\input tcilatex

\begin{document}

\begin{titlepage}


\vskip0.5truecm

\begin{center}

{\LARGE \bf A note on a standard family of twist maps}

\end{center}

\vskip  0.4truecm

\centerline {{\large Salvador Addas-Zanata}}

\vskip 0.2truecm

\centerline { {\sl Instituto de Matem\'atica e Estat\'\i stica }}
\centerline {{\sl Universidade de S\~ao Paulo}}
\centerline {{\sl Rua do Mat\~ao 1010, Cidade Universit\'aria,}} 
\centerline {{\sl 05508-090 S\~ao Paulo, SP, Brazil}}
 
\vskip 0.3truecm

\begin{abstract}
We investigate the break up of the last invariant curve for analytic families of standard maps
$$
S_\lambda :\left\{ 
\begin{array}{l}
y^{\prime }=\lambda g(x)+y \\ 
x^{\prime }=x+y^{\prime }\text{ }mod\text{ }1 
\end{array}
\right. , 
$$
where $g:S^1\rightarrow ${\rm I}\negthinspace {\rm R }is an analytic function such that $\stackunder{S^1}{\int }g(x)dx=0.$
Our main result is another evidence of how hard this
problem is. We give an example of a particular function $g$ as above such that the mapping $S_{\lambda}$ associated to
it has a "pathological" behavior.
\end{abstract} 

\vskip 0.3truecm









\vskip 3.0truecm


\hskip1.8truecm vertical rotation number, piecewise linear standard maps

\vfill
\hrule
\noindent{\footnotesize{e-mail: sazanata@ime.usp.br -  supported by FAPESP, 
grant number: 01/12449-5}}

\end{titlepage}

\section{Introduction and statement of the main result}

In this paper, we investigate the following problem:

Let $\widetilde{g}:${\rm I}\negthinspace {\rm R$\rightarrow $I}%
\negthinspace 
{\rm R} be an analytic, non-zero, periodic function, $\widetilde{g}(x+1)= 
\widetilde{g}(x)$, such that $\int_0^1\widetilde{g}(x)dx=0.$ We define the
following one parameter family ($\lambda $) of analytic diffeomorphisms of
the annulus: 
\begin{equation}
\label{slam}S_\lambda :\left\{ 
\begin{array}{l}
y^{\prime }=\lambda g(x)+y \\ 
x^{\prime }=x+y^{\prime }\text{ }mod\text{ }1 
\end{array}
\right. , 
\end{equation}
where $g:S^1\rightarrow ${\rm I}\negthinspace {\rm R }is the map induced by $
\widetilde{g}.$

For all $\lambda \in ${\rm I}\negthinspace {\rm R,} $S_\lambda $ is an
area-preserving twist mapping, because $\partial _yx^{\prime }=1,$ for any $%
(x,y)\in S^1\times ${\rm I}\negthinspace {\rm R=(}${\rm I}\negthinspace {\rm %
R/Z\negthinspace
\negthinspace Z}${\rm )}$\times {\rm I}\negthinspace {\rm R}$ and $\det
[DS_\lambda ]=1{\rm .}$ Also, the fact that $\int_0^1\widetilde{g}(x)dx=0$
implies that $S_\lambda $ is an exact mapping, which means that given any
homotopically non-trivial simple closed curve $C\subset S^1\times {\rm I}%
\negthinspace 
{\rm R,}$ the area above $C$ and below $S_\lambda (C)$ is equal the area
below $C$ and above $S_\lambda (C).$ Another obvious fact about this family
is that $S_0$ is an integrable mapping, that is, the cylinder is foliated by
invariant curves $y=y_0.$

So, KAM\ theory applies to $S_\lambda $ and we can prove that there is a
parameter $\lambda _0>0,$ such that for any $\lambda \in [0,\lambda _0]$ $%
S_\lambda $ has at least one rotational invariant curve. On the other hand,
if we choose $x_0\in S^1$ such that $g(x)\leq g(x_0)$ for all $x\in S^1,$ we
get that $S_\lambda $ does not have rotational invariant curves for all $%
\lambda \geq $$\lambda ^{*}=\frac 1{g(x_0)}>0.$ The proof of this classical
fact is very simple, so we present it here:

Given $\lambda \geq \lambda ^{*},$ choose $x_\lambda \in S^1$ such that $%
\lambda $$=\frac 1{g(x_\lambda )}.$ A computation shows that $S_\lambda
^n(x_\lambda ,0)=(x_\lambda ,n),$ for all $n\in {\rm Z\negthinspace
\negthinspace Z}.$ So there can be no rotational invariant curves.

A result due to Birkhoff implies that the set 
\begin{equation}
\label{A}A_g=\{\lambda \geq 0:\text{ }S_\lambda \text{ has at least one
rotational invariant curve}\} 
\end{equation}
is closed. So a very \ ''natural'' conjecture would be the following (see 
\cite{hall}):

\begin{description}
\item[Conjecture 1]  :$\;A_g=[0,\lambda _{cr}],$ for some $\lambda $$_{cr}>0.
$
\end{description}

Another interesting one parameter family is the following: 
\begin{equation}
\label{tlam}T_\lambda :\left\{ 
\begin{array}{l}
y^{\prime }=g(x)+y+\lambda \\ 
x^{\prime }=x+y^{\prime }\text{ }mod\text{ }1 
\end{array}
\right. 
\end{equation}

Of course $T_\lambda $ is also an area-preserving twist mapping, the
difference is that it is exact if and only if $\lambda =0,$ so when $\lambda
\neq 0$ there is no rotational invariant curve.

It can be proved (see section 2) that there is a closed interval $\rho _V=$$%
[\rho _V^{\min },\rho _V^{\max }]$ associated to $S_\lambda $ (and to $%
T_\lambda $) with the following property: Given $\omega \in \rho _V,$ there
is a point $X\in S_1\times {\rm I}\negthinspace {\rm R}$ such that%
$$
\stackunder{n\rightarrow \infty }{\lim }\frac{p_2\circ S_\lambda ^n(X)-p_2(X)%
}n=\omega , 
$$
where $p_1(x,y)=x$ and $p_2(x,y)=y$. From the exactness of $S_\lambda $ we
get that $0\in \rho _V(S_\lambda )$ for all $\lambda \in ${\rm I}%
\negthinspace {\rm R}, something that may not hold for $T_\lambda .$

In section 3 we prove a result that implies that $\rho _V^{\max }$ ($\rho
_V^{\min }$) is a continuous functions of the parameter $\lambda .$ A first
difference between $S_\lambda $ and $T_\lambda $ is that $\rho _V^{\max
}(S_\lambda )=0$ for any $\lambda \in [0,\lambda _0]$ while $\rho _V^{\max
}(T_\lambda )\neq 0$ for all $\lambda \neq 0.$ In fact, in a certain sense,
the behavior of the function $\lambda \rightarrow \rho _V^{\max }(T_\lambda
) $ is similar to the one of the rotation number of certain families of
homeomorphisms of the circle.

Given a circle homeomorphism $f:S^1\rightarrow S^1,$ a well studied family
(see for instance \cite{Herman}) is the one given by translations of $f:$ 
$$
x^{\prime }=f_\lambda (x)=f(x)+\lambda 
$$

In this case it is easy to prove that the rotation number of $f_\lambda $ is
a non-decreasing function of the parameter. We have a similar result for $%
T_\lambda $:

\begin{lemma}
\label{cresc}: $\rho _V^{\max }(T_\lambda )$ is a non-decreasing function of 
$\lambda .$
\end{lemma}

As the proof will show, this fact is an easy consequence of proposition 3,
page 466 of \cite{LeCalvez}.

If we had a similar result for $S_\lambda ,$ then conjecture 1 would
trivially be true, because $A_g=(\rho _V^{\max })^{-1}(0)$ and this set is
an interval if $\rho _V^{\max }(S_\lambda )$ is a non-decreasing function.

The main result of this note goes in the opposite direction; we present an
example in the analytic topology such that we do not know whether or not $%
A_g $ is a closed interval (although we believe it is not), but for this
example $\rho _V^{\max }(S_\lambda )$ is not a non-decreasing function of $%
\lambda $. More precisely, we have:

\begin{theorem}
\label{main}: There exists an analytic function $g^{*}$ as above such that $%
\rho _V^{\max }(S_\lambda )$ is not a non-decreasing function of $\lambda $.
\end{theorem}

The proof of the theorem implies that we can choose $g^{*}(x)=\stackunder{n=1%
}{\stackrel{N}{\sum }}a_n.\cos (2\pi nx).$ Although this choice of $g^{*}$
is a finite sum of cosines obtained as the truncation of a certain Fourier
series of a continuous function, it is still possible that for $g_S(x)=\cos
(2\pi x),$ $\rho _V^{\max }(S_\lambda )$ is in fact a non-decreasing
function, as numerical experiments suggest. Nevertheless, this shows how
subtle the problem is.

The proof of this theorem is based on a result previously obtained by the
author, on a paper due to S.Bullett \cite{bullett} on piecewise linear
standard maps and on some consequences of results from \cite{LeCalvez}.

\section{Basic tools}

First we present a theorem which is a consequence of some results from \cite
{zan}. Before we need to introduce some definitions:

1) $D_0({\rm T^2})$ is the set of torus homeomorphisms $T:{\rm %
T^2\rightarrow T^2}$ of the following form: 
\begin{equation}
\label{Ttor}T:\left\{ 
\begin{array}{l}
y^{\prime }=g(x)+y 
\text{ }mod\text{ }1 \\ x^{\prime }=x+y^{\prime }\text{ }mod\text{ }1 
\end{array}
\right. , 
\end{equation}
where $g:S^1\rightarrow {\rm I}\negthinspace {\rm R}$ is a Lipschitz
function such that $\stackunder{S^1}{\int }g(x)dx=0.$

2) $D_0(S^1\times {\rm I}\negthinspace {\rm R})$ is the set of lifts to the
cylinder of elements from $D_0({\rm T^2}),$ the same for $D_0({\rm I}%
\negthinspace {\rm R^2}).$ Given $T\in D_0({\rm T^2})$ as in (\ref{Ttor}),
its lifts $\widehat{T}\in D_0(S^1\times {\rm I}\negthinspace 
{\rm R})$ and $\widetilde{T}\in D_0({\rm I}\negthinspace 
{\rm R^2})$ write as ($\widetilde{g}$ is a lift of $g$) 
$$
\widehat{T}:\left\{ 
\begin{array}{l}
y^{\prime }=g(x)+y \\ 
x^{\prime }=x+y^{\prime }\text{ }mod\text{ }1 
\end{array}
\right. \text{ and }\widetilde{T}:\left\{ 
\begin{array}{l}
y^{\prime }= 
\widetilde{g}(x)+y \\ x^{\prime }=x+y^{\prime }\text{ } 
\end{array}
\right. 
$$

3) We say that $T\in D_0({\rm T^2})$ has a $\frac pq$-vertical periodic
orbit (set) if there is a point $A\in S^1\times {\rm I}\negthinspace {\rm R}$
such that $\widehat{T}^q(A)=A+(0,p).$ It is clear that $T^q(\pi _2(A))=\pi
_2(A),$ where $\pi _2:S^1\times {\rm I}\negthinspace {\rm R\rightarrow T^2}$
is given by $\pi _2(x,y)=(x,y$ $mod$ $1).$ The periodic orbit that contains $%
\pi _2($$A)$ is said to have vertical rotation number $\rho _V=\frac pq.$

4) Given an irrational number $\omega ,$ we say that $T\in D_0({\rm T^2})$
has a $\omega $-vertical quasi-periodic set if there is a compact $T$%
-invariant set $X_\omega \subset {\rm T^2,}$ such that for any $X\in
X_\omega $ and any $Z\in \pi _2^{-1}(X),$

$$
\rho _V(X_\omega )=\stackunder{n\rightarrow \infty }{\lim }\frac{p_2\circ 
\widehat{T}^n(Z)-p_2(Z)}n=\omega 
$$
5) We say that $T\in D_0({\rm T^2})$ has a rotational invariant curve if
there is a homotopically non-trivial simple closed curve $\gamma \subset
S^1\times {\rm I}\negthinspace {\rm R,}$ such that $\widehat{T}(\gamma
)=\gamma .$

Now we have the following:

\begin{theorem}
\label{um}: Given $T\in D_0({\rm T^2}),$ there exists a closed interval $%
0\in [\rho _V^{\min },\rho _V^{\max }]$ such that for any $\omega \in ]\rho
_V^{\min },\rho _V^{\max }[,$ there is a periodic orbit or quasi-periodic
set $X_\omega $ with $\rho _V(X_\omega )=\omega ,$ depending on whether $%
\omega $ is rational or not$.$ Moreover, $\rho _V^{\min }<0<\rho _V^{\max }$
if and only if, $T$ does not have any rotational invariant curve.
\end{theorem}

When $\omega \in \{\rho _V^{\min },\rho _V^{\max }\}$ a standard argument in
ergodic theory (see the discussion below) proves that there is an orbit with
that rotation number. In fact, much more can be said, see my forthcoming
paper \cite{meuotro}.

Following Misiurewicz and Ziemann \cite{misiu}, we can define another set
that is equal to the limit of all the convergent sequences 
$$
\left\{ \frac{p_2\circ \widehat{T}^{n_i}(Z_i)-p_2(Z_i)}{n_i},\ Z_i\in
S^1\times \negthinspace {\rm R,}\ n_i\rightarrow \infty \right\} , 
$$
which we call $\rho _V(T)^{*}.$ In the following we present a sketch of the
proof that $\rho _V(T)=\rho _V(T)^{*}.$

First note that the definition of $\rho _V(T)^{*}$ implies $\rho
_V(T)\subseteq \rho _V(T)^{*}$. Now if we define $\omega ^{-}=\inf $ $\rho
_V(T)^{*}$ and $\omega ^{+}=\sup $ $\rho _V(T)^{*},$ theorem 2.4 of \cite
{misiu} gives two ergodic $T$-invariant measures $\mu _{-}$ and $\mu _{+}$
with vertical rotation numbers $\omega ^{-}$ and $\omega ^{+},$
respectively. This means that 
$$
\stackunder{\rm T^2}{\dint }\left[ p_2\circ T(X)-p_2(X)\right] d\mu
_{-(+)}=\omega ^{-(+)}. 
$$
Therefore from the Birkhoff ergodic theorem, there are points $Z^{+}$ and $%
Z^{-}$ with $\rho _V(Z^{+})=\omega ^{+}$ and $\rho _V(Z^{-})=\omega ^{-}.$
Finally, applying theorem 6 of the appendix of \cite{zan2}, we get that $%
[\omega ^{-},\omega ^{+}]\subseteq \rho _V(T),$ so $\rho _V(T)=\rho
_V(T)^{*}.$

In the following we recall some topological results for twist maps
essentially due to Le Calvez (see \cite{lecalvez1} and \cite{lecalvez2} for
proofs), that are used in some proofs contained in this paper. Let $
\widehat{T}\in D_0(S^1\times {\rm I}\negthinspace {\rm R})$ and $\widetilde{T%
}\in D_0({\rm I}\negthinspace {\rm R^2)}$ be its lifting. For every pair $%
(s,q),$ $s\in {\rm Z\negthinspace \negthinspace Z}$ and $q\in {\rm I%
\negthinspace N^{*}}$ we define the following sets:

\begin{equation}
\label{Kpq} 
\begin{array}{c}
\widetilde{K}(s,q)=\left\{ (x,y)\in {\rm I}\negthinspace {\rm R^2}\text{: }%
p_1\circ \widetilde{T}^q(x,y)=x+s\right\} \\ \text{ and } \\ K(s,q)=\pi
_1\circ \widetilde{K}(s,q), 
\end{array}
\end{equation}
where $\pi _1:{\rm I\negthinspace R^2\rightarrow S^1\times I\negthinspace R}$
is given by $\pi _1(x,y)=(x\ mod\ 1,\ y).$

Then we have the following:

\begin{lemma}
\label{Lcal1}: For every $s\in {\rm Z\negthinspace \negthinspace Z}$ and $%
q\in {\rm I\negthinspace N^{*},\ }K(s,q)\supset C(s,q),$ a connected compact
set that separates the cylinder.
\end{lemma}

Now let us define the following functions\ on $S^1$: 
$$
\begin{array}{c}
\mu ^{-}(x)=\min \{p_2(Q) 
\text{: }Q\in K(s,q)\text{ and }p_1(Q)=x\} \\ \mu ^{+}(x)=\max \{p_2(Q)\text{%
: }Q\in K(s,q)\text{ and }p_1(Q)=x\} 
\end{array}
$$
We also have have similar functions for $\widehat{T}^q(K(s,q))$: 
$$
\begin{array}{c}
\nu ^{-}(x)=\min \{p_2(Q) 
\text{: }Q\in \widehat{T}^q\circ K(s,q)\text{ and }p_1(Q)=x\} \\ \nu
^{+}(x)=\max \{p_2(Q)\text{: }Q\in \widehat{T}^q\circ K(s,q)\text{ and }%
p_1(Q)=x\} 
\end{array}
$$

The following are important results:

\begin{lemma}
\label{graphu}: Defining Graph\{$\mu ^{\pm }$\}=\{$(x,\mu ^{\pm }(x)):x\in
S^1$\} we have: 
$$
Graph\{\mu ^{-}\}\cup Graph\{\mu ^{+}\}\subset C(s,q) 
$$

So for all $x\in S^1$ we have $(x,\mu ^{\pm }(x))\in C(s,q).$
\end{lemma}

\begin{lemma}
\label{ofpre}: $\widehat{T}^q(x,\mu ^{-}(x))=(x,\nu ^{+}(x))$ and $\widehat{T%
}^q(x,\mu ^{+}(x))=(x,\nu ^{-}(x)).$
\end{lemma}

Now we remember some ideas and results from \cite{LeCalvez}.

Given a triplet $(s,p,q)\in {\rm Z\negthinspace \negthinspace Z^2\times I%
\negthinspace N^{*},}$ if there is no point $(x,y)\in {\rm I}\negthinspace 
{\rm R^2}$ such that $\widetilde{T}^q(x,y)=(x+s,y+p),$ it can be proved that
the sets $\widehat{T}^q\circ K(s,q)$ and $K(s,q)+(0,p)$ can be separated by
the graph of a continuous function from $S^1$ to ${\rm I}\negthinspace {\rm R%
}$, essentially because from all the previous results, either one of the
following inequalities must hold: 
\begin{equation}
\label{pos}\nu ^{-}(x)-\mu ^{+}(x)>p 
\end{equation}
\begin{equation}
\label{neg}\nu ^{+}(x)-\mu ^{-}(x)<p 
\end{equation}
for all $x\in S^1,$ where $\nu ^{+},\nu ^{-},\mu ^{+},\mu ^{-}$ are
associated to $K(s,q).$

Following Le\ Calvez \cite{LeCalvez}, we say that the triplet $(s,p,q)$ is
positive (resp. negative) for $\widetilde{T}$ if $\widehat{T}^q\circ K(s,q)$
is above (\ref{pos}) (resp. below (\ref{neg})) the graph. Given $\widetilde{T%
}\in D_0({\rm I}\negthinspace {\rm R^2}),$ we have:

$$
\widetilde{T}(x,y)=(x^{\prime },y^{\prime })\Leftrightarrow y=m(x,x^{\prime
})\text{ and }y^{\prime }=m^{\prime }(x,x^{\prime }), 
$$
where $m$ and $m^{\prime }$ are continuous maps from ${\rm I}\negthinspace 
{\rm R^2}$ to ${\rm I}\negthinspace {\rm R}$ with some especial properties$%
{\rm .}$

If $\widetilde{T},\widetilde{T^{*}}\in D_0({\rm I}\negthinspace {\rm R^2}),$
we say that $\widetilde{T}\leq \widetilde{T^{*}}$ if $m^{*}\leq m$ and $%
m^{\prime }\leq m^{*\prime },$ where $(m,m^{\prime })$ is associated to $
\widetilde{T}$ and $(m^{*},m^{*\prime })$ to $\widetilde{T^{*}}.$

\begin{proposition}
\label{compte}: If $(s,p,q)$ is a positive (resp. negative) triplet of $
\widetilde{T}$ and if $\widetilde{T}\leq \widetilde{T^{*}}$ (resp. $
\widetilde{T}\geq \widetilde{T^{*}}$), then $(s,p,q)$ is a positive (resp.
negative) triplet of $\widetilde{T^{*}}.$
\end{proposition}

Now we present an amazing example of a twist homeomorphism from $D_0({\rm T^2%
}).$ First, let $g^{\prime }:S^1\rightarrow {\rm I}\negthinspace {\rm R}$ be
given by $g^{\prime }(x)=\left| x-\frac 12\right| -\frac 14$ and so the lift 
$\widetilde{g}^{\prime }:{\rm I}\negthinspace {\rm R\rightarrow I}%
\negthinspace {\rm R}$ is continuous, $\widetilde{g}^{\prime }(x+1)= 
\widetilde{g}^{\prime }(x),$ $\int_0^1\widetilde{g}^{\prime }(x)dx=0,$ $Lip( 
\widetilde{g}^{\prime })=1$ and $\widetilde{g}^{\prime }(x)=\widetilde{g}%
^{\prime }(-x).$ Also, $\widetilde{g}^{\prime }$ is differentiable
everywhere, except at points of the form $\frac n2,$ $n\in {\rm Z%
\negthinspace \negthinspace Z.}$ The one parameter family $S_\lambda
^{\prime }\in D_0({\rm T^2})$ is given by: 
\begin{equation}
\label{slaml}S_\lambda ^{\prime }:\left\{ 
\begin{array}{l}
y^{\prime }=\lambda g^{\prime }(x)+y 
\text{ }mod\text{ }1 \\ x^{\prime }=x+y^{\prime }\text{ }mod\text{ }1 
\end{array}
\right. 
\end{equation}

In \cite{bullett} this family is studied in detail and among other things,
the following theorem is proved:

\begin{theorem}
\label{bull}: There are no rotational invariant curves for $S_\lambda
^{\prime }$ when

$\lambda \in ]0,918,1[\cup ]4/3,\infty [$ and for $\lambda =4/3$ there are
''lots '' of rotational invariant curves.
\end{theorem}

\section{Proofs}

\subsection{Preliminary results}


{\it Proof of lemma \ref{cresc}:}

\vskip0.2truecm

This result is a trivial consequence of proposition \ref{compte}. Given $%
\lambda _1<\lambda _2,$ we get from expression (\ref{tlam}) that $\widetilde{%
T}_{\lambda _1}\leq \widetilde{T}_{\lambda _2}.$ So if $\rho _V^{\max
}(T_{\lambda _2})<p/q<\rho _V^{\max }(T_{\lambda _1})$ for a certain
rational number $p/q,$ then for any $s\in {\rm Z\negthinspace 
\negthinspace Z}$ the triplet $(s,p,q)$ is negative for $\widetilde{T}%
_{\lambda _2}$, which implies by proposition \ref{compte} that it is also
negative for $\widetilde{T}_{\lambda _1},$ which contradicts the fact that $%
\rho _V^{\max }(T_{\lambda _1})>p/q.$ 

\vskip0.5truecm

Now we prove the following theorem that has its own interest. It is easy to
see from the proof that it is valid in a more general context.

\begin{theorem}
\label{MZ}:\ The functions $\rho _V^{\max },\rho _V^{\min }:D_0({\rm T^2}%
)\rightarrow {\rm I}\negthinspace {\rm R}$ are continuous.
\end{theorem}

Remark: The proofs are analogous, so we do it only for $\rho _V^{\max }.$

\vskip0.2truecm


{\it Proof:}

Suppose that there is a $T_0\in D_0({\rm T^2})$ such that $\rho _V^{\max }$
is not continuous at $T_0.$ This means that there is an $\epsilon >0$ and a
sequence $D_0({\rm T^2})\ni $$T^n\stackrel{n\rightarrow \infty }{\rightarrow 
}T_0$ in the $C^0$  topology, such that either:

1) $\rho _V^{\max }(T^n)>\rho _V^{\max }(T_0)+\epsilon ,$ for all $n$, or

2) $\rho _V^{\max }(T^n)<\rho _V^{\max }(T_0)-\epsilon ,$ for all $n.$

The first possibility means that there exists a rational number $p/q$ such
that $\rho _V^{\max }(T^n)>p/q>\rho _V^{\max }(T_0).$ This implies that for
any $s\in {\rm Z\negthinspace \negthinspace Z,}$ the triplet $(s,p,q)$ is
non-negative for $\widetilde{T}^n$ (as the value of $s$ is irrelevant in
this setting, we fix $s=0$). But as $\rho _V^{\max }(T_0)<p/q,$ $(0,p,q)$ is
negative for $\widetilde{T}_0.$ As $T^n\stackrel{n\rightarrow \infty }{%
\rightarrow }T_0,$ we get from the upper semi-continuity in the Hausdorff
topology of the maps  
\begin{equation}
\label{hauscon}T\rightarrow K(0,q)\text{ and }T\rightarrow \widehat{T}%
^q(K(0,q)) 
\end{equation}
that $(0,p,q)$ is a negative triplet for all mappings sufficiently close to $
\widetilde{T}_0,$ which is a contradiction.

In the same way, the second possibility means that there exists a rational
number $p/q$ such that $\rho _V^{\max }(T^n)<p/q<\rho _V^{\max }(T_0).$ This
implies that there exists $Q\in C(0,q)$ such that 
\begin{equation}
\label{certo}p_2\circ \widehat{T_0}^q(Q)-p_2(Q)>p. 
\end{equation}

Now we prove the following claim, which implies the theorem:

\begin{description}
\item[Claim]  :\ Any mapping $T\in D_0({\rm T^2})$ sufficiently close to $T_0
$ will satisfy an inequality similar to (\ref{certo}).
\end{description}


{\it Proof:}

First of all, let us define $P_0=(x_Q,\mu ^{-}(x_Q)),$ where $x_Q=p_1(Q).$
From lemma \ref{ofpre} and the definition of $\mu
^{-} $ and $\nu ^{+},$ we get that $\nu ^{+}(x_Q)=p_2\circ \widehat{T_0}%
^q(P_0)>p_2(P_0)+p=\mu ^{-}(x_Q)+p.$ So there exists $\delta >0$ such that
for any $Z\in \overline{B_\delta (P_0)}$ we have 
$$
p_2\circ \widehat{T_0}^q(Z)>p_2(Z)+p. 
$$
Therefore, there exists a neighborhood $T_0\in {\cal U}\subset D_0({\rm T^2}%
) $ in the $C^0$ topology such that for any $T\in {\cal U},$ we get $%
p_2\circ \widehat{T}^q(Z)>p_2(Z)+p,$ for all $Z\in B_\delta (P_0).$ Now
defining $\overline{AB}=\left\{ x_Q\times {\rm I}\negthinspace {\rm R}%
\right\} {\rm \cap }B_\delta (P_0),$ lemma \ref{graphu} implies that if we
choose a sufficiently small neighborhood $V$ of $C(0,q),$ then for all
homotopically non-trivial simple closed curves $\gamma \subset V,$ we get
that $\gamma \cap \overline{AB}\neq \emptyset .$ By the upper
semi-continuity in the Hausdorff topology of the maps in (\ref{hauscon}), if
we choose a sufficiently small sub-neighborhood ${\cal U^{\prime }}{\it %
\subset }{\cal U}$ we get for any $T\in {\cal U^{\prime }}$ that the set $%
C(0,q)$ associated to $T$ is also contained in $V.$ Therefore it must cross $
\overline{AB}.$

So given any mapping $T\in {\cal U^{\prime }}{\it \subset }{\cal U}$, there
is a point $Q^{\prime }\in C(0,q)\cap \overline{AB}$ which therefore
satisfies $p_2\circ \widehat{T}^q(Q^{\prime })>p_2(Q^{\prime })+p.$ 
\vskip0.3truecm

Finally, the above claim implies that $\rho _V^{\max }(T^n)\geq p/q$ for
sufficiently large $n$, which is a contradiction. 

\subsection{Main theorem}

In this section we prove theorem (\ref{main}).

\vskip0.2truecm

First of all we note that from theorem (\ref{bull}), the mapping $S_\lambda
^{\prime }\in D_0({\rm T^2})$ (see (\ref{slaml})) has no rotational
invariant curve for $\lambda =0.95$ and has ''lots'' of rotational invariant
curves for $\lambda =4/3$ . Using theorem (\ref{um}) one gets that $\rho
_V^{\max }(S_{0.95}^{\prime })=\epsilon >0$ and $\rho _V^{\max
}(S_{4/3}^{\prime })=0.$ A classical result in Fourier analysis implies that
the Fourier series $\widetilde{g}_N^{\prime }(x)=\stackrel{N}{\stackunder{n=1%
}{\sum }}a_n\cos (2\pi nx)$ of $\widetilde{g}^{\prime }$ converges uniformly
to $\widetilde{g}^{\prime }$. So if we choose $N>0$ sufficiently large, we
get from theorem (\ref{MZ}) that $\rho _V^{\max }(S_{N,0.95}^{\prime
})>\epsilon /2$ and $\rho _V^{\max }(S_{N,4/3}^{\prime })<\epsilon /10,$
where $S_{N,\lambda }^{\prime }$ is the twist mapping associated to $%
g_N^{\prime }$. 

\end{document}